\documentclass{conm-p-l}
\usepackage{graphicx}

\newtheorem{theorem}{Theorem}[section]

\newtheorem{conjecture}[theorem]{Conjecture}
\newtheorem{proposition}[theorem]{Proposition}

\theoremstyle{definition}

\theoremstyle{remark}

\numberwithin{equation}{section}

\newcommand{\abs}[1]{\lvert#1\rvert}

\newcommand\R{{\mathbb R}}
\renewcommand\O{{\mathcal O}}

\newcommand\lfs{\operatorname{lfs}}

\newcommand\half{\tfrac12}

\let\paragraph=\S
\renewcommand\S{{\mathbb S}}

\def\incgrh#1#2{\includegraphics[height=#2in]{#1}}
\def\incgrs#1#2{\includegraphics[scale=#2]{#1}}
\def\endfig#1#2#3{\caption[#2]{#3}\label{fig:#1}\end{figure}}

\def\figs#1#2#3#4%
{ \begin{figure}
  \begin{center} \incgrs{#1}{#2} \end{center}
  \endfig{#1}{#3}{#4} }

\def\figh#1#2#3#4%
{ \begin{figure}
  \begin{center} \incgrh{#1}{#2} \end{center}
  \endfig{#1}{#3}{#4} }

\def\figtwo#1#2#3#4%
{ \begin{figure}
  \begin{center} \incgrh{#1A}{#2} \hspace{.2in} \incgrh{#1B}{#2} \end{center}
  \endfig{#1}{#3}{#4} }

\def\figr#1{Figure~\ref{fig:#1}}

\newcommand{\Rp}{R^p}

\begin{document}

\title{Approximating Ropelength by Energy Functions}

\author{John M. Sullivan}
\address{Department of Mathematics, University of Illinois,
Urbana, IL 61801--2975}
\email{jms@math.uiuc.edu}
\urladdr{http://www.math.uiuc.edu/\textasciitilde{}jms/}
\thanks{The author was supported in part by NSF grant DMS-00-71520,
and would like to thank Jason Cantarella, Rob Kusner and Eric Rawdon
for useful conversations, and the anonymous referee for many
helpful suggestions on the first version of this paper.}
\copyrightinfo{2002}{John M. Sullivan}

\date{\today.}

\begin{abstract}
 The ropelength of a knot is the quotient of its length by its
 thickness.  We consider a family of energy functions $\Rp$ for knots,
 depending on a power $p$, which approach ropelength as $p$ increases.
 We describe a numerically computed trefoil knot which seems to
 be a local minimum for ropelength; there are nearby critical points
 for $\Rp$, which are evidently local minima for large enough $p$.
\end{abstract}

\maketitle

\section{Thickness and ropelength}

We measure the \emph{ropelength} of a knot as the quotient of its length
by its \emph{thickness}.
The thickness is the radius of the largest embedded normal tube
(called the thick tube) around the knot.
Ropelength is a mathematical model of how much physical rope it would
take to tie the knot.  In~\cite{cks2} we showed
that in any knot or link type there is a ropelength minimizer,
and that minimizers are necessarily $C^{1,1}$ curves.
But still, explicit examples of tight (ropelength minimizing) links
are known only in very special cases.

One way to define the thickness of a link (following~\cite{GM})
is to consider circles through three points on the link.
For any three distinct points $x$, $y$, $z$ in $\R^3$, we let $r(x,y,z)$ be
the radius of the (unique) circle through these points (setting
$r=\infty$ if the points are collinear).  Also, if $T_x$ is a unit
vector at~$x$, we let $r(T_x,y)$ be the radius of the circle through~$y$
tangent to~$T_x$ at~$x$.  Note that
$$r(T_x,y)=\frac{\|x-y\|^2}{2\sqrt{\|x-y\|^2-(T_x\cdot(x-y))^2}}.$$

Now let $L$ be a link in~$\R^3$, that is, a disjoint union of simple
closed curves.  We define the \emph{thickness} $\tau(L)$ of~$L$ as
$$\tau(L):=\!\!\inf_{\substack{x,y,z\in L\\x\ne y\ne z\ne x}} \!\!r(x,y,z).$$
This quantity vanishes when $L$ fails to be $C^{1,1}$, so from now
on we will restrict attention to $C^{1,1}$ links.

The infimum in the definition of~$\tau$ is always achieved when
(at least) two of the three points approach each other~\cite{GM,cks2}: we have
$$\tau(L)=\inf_{y\ne x} r(T_x L,y),$$
where $T_xL$ denotes the tangent vector to the $C^1$ curve $L$ at $x$.
We can equivalently~\cite{cks2} define
the thickness in terms of the \emph{medial axis} (see~\cite{ACK}) of~$L$.
The medial axis is the set of points in space which fail to have a unique
nearest point in~$L$.  Then the thickness $\tau(L)$ is exactly
Federer's \emph{reach}~\cite{Federer} of~$L$,
the distance from~$L$ to its medial axis.

Gonzalez and Maddocks~\cite{GM} considered the quantity
$$\tau^{\mathrm{GM}}_x(L) := \inf_{y,z} r(x,y,z) = \inf_y r(T_y L,x).$$
which they called the ``global radius of curvature'' but which could well be
termed the \emph{local thickness} of the link~$L$ at the point~$x\in L$.
However, it seems more natural to define local thickness in terms of
higher-order contact at~$x$:
$$\tau_x(L) := \inf_{y\ne x} r(T_x L, y).$$
Either definition is a local thickness in the sense that
$$\tau(L)=\inf_x\tau_x(L)=\inf_x\tau^{\mathrm{GM}}_x(L).$$
However, $\tau_x(L)$ has a nice geometric interpretation as the
radius of the largest sphere which, staying tangent to $L$ at $x$,
can be rotated completely around $L$ without touching $L$ at any other point.

Given a circle~$C$ of radius~$r$, consider the set of points in space
whose distance to~$C$ is less than~$r$.
Following Dan Asimov, we will call this a \emph{bialy}
(or more precisely an open solid bialy),
with \emph{neck} at the center of~$C$ and \emph{axis} normal to~$C$.
A bialy is an open solid torus of revolution whose major and minor
radii are equal.  Given a point~$x$ on~$L$, the bialys with neck at~$x$
and axis along $T_xL$ are nested, and it is clear that
$\tau_x(L)$ is radius of the largest one which avoids $L$.

We now check that our definition of local thickness also agrees
with the notion of \emph{local feature size} from the theory of medial axes.
Let $r(T_xL, T_yL)$ denote the radius of the sphere
tangent to both $T_xL$ and~$T_yL$
(the smallest such sphere if these vectors are cocircular).
Then it is clear from our discussion above that
$$\tau_x(L) = \inf_y \, r(T_xL, T_yL),$$
since the interiors of all spheres of radius~$r$ tangent at~$x$
fill out the bialy at~$x$.

A formal definition of medial axis is the following.  Suppose $L$
is a compact subset of~$\R^3$.  Then for any point $p\in\R^3$, we can consider
the distance to~$p$ as a function $d_p:L\to\R$.  It achieves its
minimum $d(p,L)$ on some set $\min(d_p)\subset L$.  Then the medial
axis is $M(L):=\{p\in\R^3: \#\min(d_p) > 1\}$, the set of~$p$ for
which this minimum is achieved at more than one point.  At a point
$p\in M$, the \emph{local feature size} is simply the
distance $d(p,L)$ to~$L$.
At a point $x\in L$, the \emph{local feature size} is
$$\lfs(x) := \inf\{d(p,L): p\in M, x\in \min(d_p)\}.$$
(We note that this definition, intended for use when $L$ is a curve
or surface, is only one of many inequivalent definitions for
local feature size that have been given in the computer science literature.)

In terms of these definitions, we are now ready to
state our first proposition:

\begin{proposition}
The local thickness $\tau_x(L)$ of a link $L$ is also the
local feature size of~$L$ at~$x$.
\end{proposition}
\begin{proof}
The medial axis $M(L)$ is exactly the set of centers $p$ of spheres
which are tangent to~$L$ at two or more points, and whose interiors
avoid $L$.  The local feature size at $p\in M$ (or at the points of
tangency along~$L$) is the radius of the corresponding sphere.
But the local thickness at~$x$ is also the radius of the smallest
bitangent sphere tangent at~$x$.
\end{proof}

\section{Smooth approximations to ropelength}

We can approximate the infimum in the definition of ropelength
with $L^p$ energies.  For a link $L$ of length $1$, we define
$$\Rp_x(L) := \int_{y\in L} r(T_x L, y)^{-p}\,ds,
\qquad \Rp(L) := \left(\int_x \Rp_x(L)\,ds\right)^{\!1/p}.$$
We extend this energy to links of arbitrary length so
that it is scale-invariant.  Aside from this rescaling,
the energy $\Rp(L)$ is the same as the energy $U_{p,2}(L)$
suggested in~\cite{GM}.

Then clearly $\Rp(L)$ is bounded above by the ropelength $R(L)$,
and in fact $$\lim_{p\to\infty} \Rp(L) = R(L).$$
For $p\ge 1$, it also follows from H\"older's inequality that
$\Rp(L) \ge R^1(L)$.
The examples in~\cite{cks2} show that tight links
(ropelength minimizers) need not be smoother than $C^{1,1}$.
However, since the energies $\Rp$ smooth out
the hard-shell potential involved in a thickness constraint,
we expect that minimizers for the energies $\Rp$ should be smooth.

Diao, Ernst and Janse van Rensburg~\cite{DEJvR-ideal}
have listed several properties that one might desire for
knot-energy functionals.  For any $p\ge1$, we conjecture that
the energy $\Rp$ is ``basic'' in the sense that the absolute minimum
is achieved uniquely by the round circle (for which $\Rp=2\pi$).
It would suffice to prove this conjecture for the special case $p=1$.

For $p<2$, the quantity $\Rp$ is not very useful as a knot energy,
because it does not have an infinite barrier to self-crossings.
It is straightforward to check that the energy of two skew lines
which approach a right-angled crossing becomes infinite if and only
if $p\ge 2$.

The following conjecture would imply that, for $p\ge 2$,
$\Rp$ is ``charge'' and ``tight'', meaning that it approaches infinity
for sequences of links which approach a curve with a self-intersection,
or in which a knotted arc shrinks to zero size.
\begin{conjecture}
For $p\ge 2$, the energy $\Rp$ is bounded below by some monotonic
function of ropelength (so $\Rp$ approaches infinity for any
sequence of links with fixed length and thickness approaching zero).
\end{conjecture}
We will prove part of this conjecture, enough to conclude
that $\Rp$ is ``charge''.

\begin{proposition}
If $L_0$ is the smooth limit of smooth links $L_n$,
and $L_0$ has a self-crossing,
then for $p>2$ the energies $\Rp(L_n)$ approach infinity.
\end{proposition}
\begin{proof}
Let $K/2$ be an upper bound for the curvature of~$L_0$,
so that $K$ bounds the curvature of~$L_n$ for large $n$.
There is also a uniform bound $T$ on the third derivatives.
The thickness of~$L_n$ approaches $0$.

Let $L$ be any of the $L_n$ with $n$ large enough that
its thickness $\tau$ satisfies
$$\tau<\frac1{4K},\qquad \tau << 1/T.$$
Since $\tau<1/K$, the thickness is achieved by a doubly critical pair of
points $x,y$, with $L$ perpendicular to the segment $\overline{xy}$
at both ends.  Picking an appropriate coordinate system,
we can assume that $x=(0,0,\tau)$ and $y=(0,0,-\tau)$,
with $T_xL=(1,0,0)$ and $T_yL=(\cos\phi,\sin\phi,0)$ for some angle $\phi$.
We can then expand $L$ in Taylor series near $x$ and $y$ as:
\begin{eqnarray*}
\alpha(s) &=& (s, k_1 s^2, \tau+k_2 s^2)+\O(s^3),\\
\beta(t) &=& (t\cos\phi + k_3 t^2\sin\phi,
            t\sin\phi - k_3 t^2\cos\phi,
           -\tau - k_4 t^2)+\O(t^3),
\end{eqnarray*}
where the $k_i$ are components of the curvature vector, so $\abs{k_i}<K$,
and the constants in the first omitted terms depend only on $T$.
The tangent vector at $x=\alpha(s)$ is $(1, 2k_1s, 2k_2s)+\O(s^2)$.
Thus for $y=\beta(t)$ we can compute
$$r(T_xL,y) =
  \tau + s^2\big(\tfrac1{4\tau}-k_2\big)
       + t^2\big(\tfrac1{8\tau}(2-\sin^2\phi)+\tfrac{k_4}2\big)
       - st\big(\tfrac{\cos\phi}{2\tau}+k_2\big) +\O(s^3,t^3).$$
This gives
$$r(T_xL,y) \le \tau \left(1 + C(\abs{s}+\abs{t})^2 +\O(s^3,t^3)\right),$$
where $C= 1/(4\tau^2)+ K/\tau$.  Since $K\tau<1/4$, we have $C<1/(2\tau^2)$.
Now pick some small length $a=\tau^q$,
and consider the integral defining $\Rp(L)$
taken over just the ranges $0\le s,t \le a$.  We get
\begin{eqnarray*}
\Rp(L)^p &>& \int_0^a\!\!\int_0^a \!r^{-p}\,ds\,dt \\
 &=& \tau^{-p} \iint \!\left(1-pC(s+t)^2+\O(s^3,t^3)\right)\,ds\,dt \\
 &=& \tau^{-p} \left(a^2 - pC(\tfrac76 a^4) + \O(a^5)\right) \\
 &\ge& \tau^{-p} \left(\tau^{2q} - p\tau^{4q}/\tau^2 + \O(\tau^{5q}) \right)
\end{eqnarray*}
Now if $q>1$, the first term dominates the others.
(We should note that a more careful analysis would confirm
that the constants in the term of order $\tau^{5q}$ are still
uniformly bounded in terms of~$T$.)
For $1<q<p/2$, this expression diverges to infinity as $\tau\to0$, as desired.
\end{proof}

For $p\ge 2$, we expect that $\Rp$ is also ``strong'',
in the sense of allowing only a finite number of
link types under any given energy level, but this seems hard to prove.

\section{Numerical simulations}

We have earlier reported~\cite{KS-knot} on numerical simulations,
conducted using Brakke's Evolver~\cite{evolver},
of knots minimizing the M\"obius-invariant energy of~\cite{FHW}.
We have recently used the Evolver for
numerical simulations of tight (ropelength minimizing) links,
using a discretization, for polygonal approximations,
of the energies $\Rp$ introduced here.
These simulations will be described in more detail in a forthcoming paper.

As an example of the results, we describe here a trefoil
knot, shown in \figr{t32},
\figtwo{t32}{2.3}{Projections of the local minimum trefoil.} 
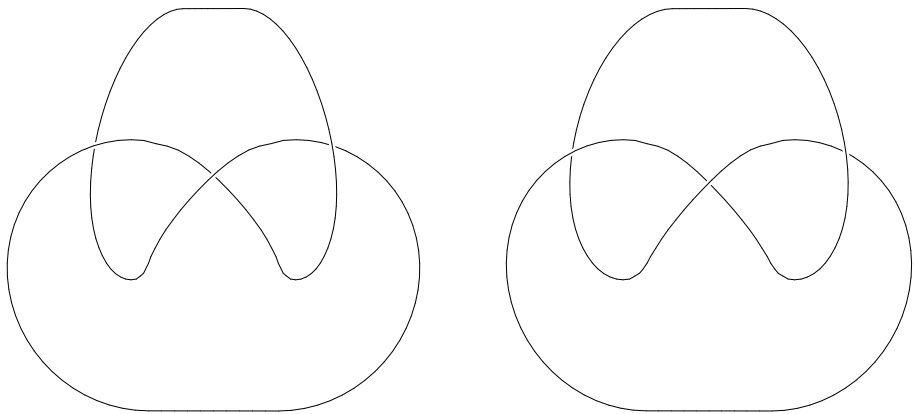
{\small A pair of orthogonal projections
(printed for cross-eyed stereo viewing)
of the trefoil knot $\gamma_p$ for $p=4096$.
This is a numerically computed local minimum for the
energy $\Rp$, and it is presumably close to a local minimum
for ropelength.  Its ropelength is about $37.5$, as opposed
to $32.7$ for the minimizing trefoil.}
which seems to be a local mininum for ropelength.
This trefoil is embedded in space as a $(3,2)$--torus knot,
as opposed to the $(2,3)$--torus knot which gives the (presumed) global minimum.
Its symmetry group in $\R^3$ is $2\,2\,2$ in the Conway/Thurston
notation for $2$-orbifolds: that is, it is generated by the rotations of
order $2$ around three orthogonal axes.

Numerically we have computed a family of very similar curves $\gamma_p$
with this symmetry which are critical points for the energies $\Rp$.
Presumably, they are global minima if the symmetry is fixed.  But without
enforcing the symmetry, they seem to be unstable saddle points for
$\Rp$ when $p\le512$ but stable local minima when $p\ge 1024$.
(We have only tested powers $p$ which are powers of~$2$.)
In \figr{t32-kt} we plot the curvature and torsion of $\gamma_{4096}$,
as functions of arclength along the curve.

Presumably, the limit of these curves as $p\to\infty$ is a curve $\gamma$
which is a local minimum for ropelength.  This limit curve is interesting
in that it seems to include two straight segments.  The thick tube around
these arcs does not contact that around any other part of~$\gamma$.
Presumably, at the ends of these segments, where the tubes do contact,
the curvature of~$\gamma$ jumps discontinously to some positive value
(close to $1$ if we normalize to $\tau=1$).
Thus $\gamma$ would be an example of a ropelength-critical knot which is not
$C^2$ smooth.

The only links so far proved to be tight are certain examples from~\cite{cks2}.
In these links, each component is planar, is piecewise circular (or straight),
but fails to be $C^2$ unless it is a circle.  No explicit tight knots
are known.  It seems unlikely that any of them would be more
than piecewise smooth.  But our conjectured local-minimum trefoil
gives the clearest example yet where a discontinuity in curvature
should happen.
\figtwo{t32-kt}{2.3}{Curvature and torsion of the local minimum trefoil.}
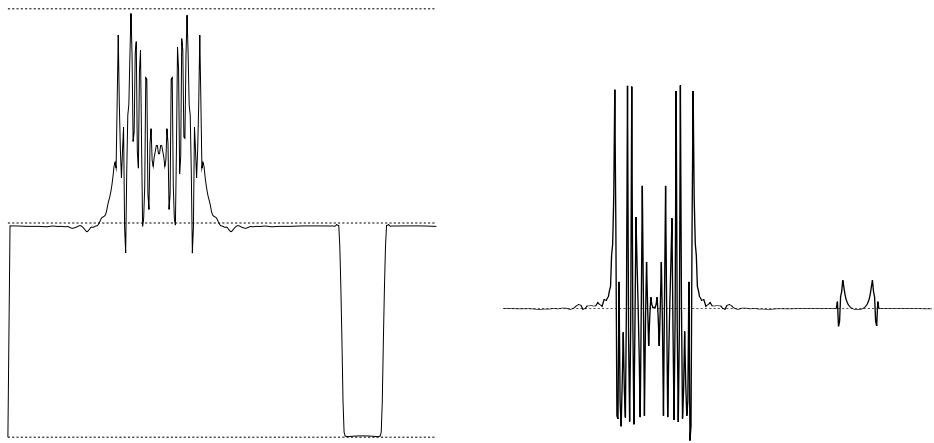
{\small Numerically computed curvature (left) and torsion (right) functions
of the curve $\gamma_{4096}$, plotted against arclength,
going slightly more than halfway around the curve.
The high frequencies should probably be viewed as numerical noise.
There is presumably a similar $(3,2)$--torus knot $\gamma$ which is
a local minimum for ropelength.
The curvature of~$\gamma$ seems to jump abruptly from $0$ to $\half$,
and then smoothly increase further.
The torsion is close to zero except where the curvature is high.
Note that for the finite value $p=4096$ shown, the arcs which
would be straight in $\gamma$ are not exactly straight.}

\bibliographystyle{hamsalpha} \bibliography{thick}
\end{document}